\documentclass[a4paper, 12pt]{article}%
\usepackage{amsfonts}
\usepackage{fancyhdr}
\usepackage[a4paper, top=3.0cm, bottom=3.0cm, left=2.0cm, right=2.0cm]%
{geometry}
\usepackage{amssymb}
\usepackage{graphicx}
\usepackage{amsmath}%

\newtheorem{theorem}{Theorem}

\newtheorem{corollary}[theorem]{Corollary}

\newtheorem{lemma}[theorem]{Lemma}

\newtheorem{proposition}[theorem]{Proposition}

\begin{document}

\title{A Spectral Strong Approximation Theorem for Measure Preserving Actions}
\author{Mikl\'{o}s Ab\'{e}rt }
\maketitle

\begin{abstract}
Let $\Gamma$ be a finitely generated group acting by probability measure
preserving maps on the standard Borel space $(X,\mu)$. We show that if
$H\leq\Gamma$ is a subgroup with relative spectral radius greater than the
global spectral radius of the action, then $H$ acts with finitely many ergodic
components and spectral gap on $(X,\mu)$. This answers a question of Shalom
who proved this for normal subgroups.

\end{abstract}

\section{Introduction}

Let $(X,\mathcal{B},\mu)$ be a standard Borel probability space. Let
$L^{2}(X)=L^{2}(X,\mathcal{B},\mu)$ denote the space of square integrable
measurable real functions on $X$ and let $L_{0}^{2}(X)\subseteq L^{2}(X)$ be
the subspace of functions with zero integral. Let $\mathrm{Aut}(X,\mu)$ denote
the Polish group of $\mu$-preserving Borel isomorphisms of $(X,\mathcal{B}%
,\mu)$.

Let $\lambda$ be a Borel probability measure on $\mathrm{Aut}(X,\mu)$. One can
associate the averaging operator $M_{\lambda}:L^{2}(X)\rightarrow L^{2}(X)$
defined by%

\[
(fM_{\lambda})(x)=%
{\displaystyle\int\limits_{a\in\mathrm{Aut}(X,\mu)}}
f(xa)d\lambda(a)\text{ \ (}f\in L^{2}(X),x\in X\text{)}%
\]
Let the top of the spectrum of $\lambda$ be
\[
\rho^{+}(\lambda)=\rho^{+}(X,\mathcal{B},\mu,\lambda)=\sup_{0\neq f\in
L_{0}^{2}(X)}\frac{\left\langle fM_{\lambda},f\right\rangle }{\left\langle
f,f\right\rangle }%
\]
and let the norm of $\lambda$ be
\[
\rho(\lambda)=\left\Vert M_{\lambda}\right\Vert =\sup_{0\neq f\in L_{0}%
^{2}(X)}\sqrt{\frac{\left\langle fM_{\lambda},fM_{\lambda}\right\rangle
}{\left\langle f,f\right\rangle }}\text{.}%
\]
Then we have
\[
0\leq\rho^{+}(\lambda)\leq\rho(\lambda)\leq1\text{.}%
\]
The operator $M_{\lambda}$ is self-adjoint when $\lambda$ is \emph{symmetric},
that is, $\lambda=\lambda^{-1}$ where $\lambda^{-1}$ is obtained by composing
$\lambda$ with the inverse map. In general, we have $\rho(\lambda)^{2}%
=\rho^{+}(\lambda\lambda^{-1})$. The inequality $\rho^{+}(\lambda)\geq0$ is
proved in Proposition \ref{nulla}.

Let $\Gamma$ be a countable group and let $\lambda$ be a symmetric probability
measure on $\Gamma$ such that the support of $\lambda$ generates $\Gamma$. The
averaging operator $M_{\lambda}$ now naturally acts on $l^{2}(\Gamma)$ and is
self-adjoint. Let the \emph{spectral radius of }$\lambda$ be
\[
\rho(\lambda)=\sup_{0\neq f\in l^{2}(\Gamma)}\frac{\left\langle fM_{\lambda
},f\right\rangle }{\left\langle f,f\right\rangle }=\left\Vert M_{\lambda
}\right\Vert =\lim_{n\rightarrow\infty}\sqrt[2n]{p_{e,e,2n}}%
\]
where $p_{e,e,n}$ is the probability of return for the $\lambda$-random walk
on $\Gamma$ in $n$ steps. 

A \emph{p.m.p. action }$\varphi$\emph{ of }$\Gamma$\emph{\ on }$(X,B,\mu)$ is
a homomorphism from $\Gamma$ to $\mathrm{Aut}(X,\mu)$. The push-forward
$\varphi(\lambda)$ is a Borel probability measure on $\mathrm{Aut}(X,\mu)$ and
the action $\varphi$ has spectral gap if and only if $\rho^{+}(\varphi
(\lambda))<1$. In Proposition \ref{locglob} we show that
\[
\rho(\lambda)\leq\rho^{+}(\varphi(\lambda))\text{.}%
\]

Let $H$ be an arbitrary subgroup of $\Gamma$ and let $\lambda$ be a symmetric
probability measure on $\Gamma$. Then $M_{\lambda}$ also acts on
$l^{2}(H\backslash\Gamma)$ where $H\backslash\Gamma$ is the set of right
cosets of $H$ in $\Gamma$. Let the \emph{relative spectral radius}
\[
\rho(\Gamma,H,\lambda)=\sup_{0\neq f\in l^{2}(H\backslash\Gamma)}%
\frac{\left\langle fM_{\lambda},f\right\rangle }{\left\langle f,f\right\rangle
}=\lim_{n\rightarrow\infty}\sqrt[2n]{p_{e,H,2n}}=\left\Vert M_{\lambda
}\right\Vert
\]
where $p_{e,H,n}$ is the probability that the $\lambda$-random walk of length
$n$ starting at $e$ ends in $H$. We have
\[
\rho(\lambda)\leq\rho(\Gamma,H,\lambda)\leq1
\]
The quantity $\rho(\Gamma,H,\lambda)$ can be thought of as the \emph{dimension
of }$H$\emph{ in }$\Gamma$.

Our first theorem says that when $H\leq\Gamma$ is too big compared to a p.m.p.
action $\varphi$ of $\Gamma$, it can not effectively `hide in the action'. For
a finite symmetric set $S$ of $\Gamma$ let $\lambda_{S}$ denote the uniform
probability measure on $S$.

\begin{theorem}
\label{fotetel}Let $\Gamma$ be generated by the finite symmetric set $S$ and
let $\varphi$ be a p.m.p. action of $\Gamma$. Then for every subgroup $H$ of
$\Gamma$ with
\[
\rho(\Gamma,H,\lambda_{S})>\rho^{+}(\varphi(\lambda_{S}))
\]
there exists a finitely generated subgroup $H^{\prime}$ of $H$ such that
$H^{\prime}$ has finitely many ergodic components and $H^{\prime}$ acts on
each component with spectral gap.
\end{theorem}

Theorem \ref{fotetel} answers a question of Shalom, who proved it in the case
when $H$ is normal in $\Gamma$ \cite{shalom}. We use a different approach from
Shalom, which also provides an explicit generating set for $H^{\prime}$ and
effective bounds on the expansion properties of this generating set.

The proof of Theorem \ref{fotetel} uses geometric expansion. Let $\lambda$ be
a Borel probability measure on $\mathrm{Aut}(X,\mu)$. For a Borel subset
$Y\in\mathcal{B}$ let
\[
e_{\lambda}(Y)=%
{\displaystyle\int\limits_{a\in\mathrm{Aut}(X,\mu)}}
\mu(Ya\diagdown Y)d\lambda(a)=\left\langle \chi_{Y}M_{\lambda},\chi_{Y^{c}%
}\right\rangle \text{.}%
\]
be the probability that a $\lambda$-random edge starting at $Y$ leaves $Y$.
Let the \emph{expansion constant of }$\lambda$ be
\[
\mathrm{h}(\lambda)=\mathrm{h}(X,\mathcal{B},\mu,\lambda)=\inf\left\{
\frac{e_{\lambda}(Y)}{\mu(Y)(1-\mu(Y))}\mid Y\in\mathcal{B}\text{, }%
0<\mu(Y)<1\right\}  \text{.}%
\]
We call the system $\lambda$ an \emph{expander} if $\mathrm{h}(\lambda)>0$.
Adapting Cheeger's inequalities for group actions in \cite{lynaz} gives the
estimates
\[
1-\mathrm{h}(\lambda)\leq\rho^{+}(\lambda)\leq1-\mathrm{h}(\lambda
)^{2}/8\text{.}%
\]
In particular, $\lambda$ is an expander if and only if $\rho^{+}(\lambda)<1$.
See Proposition \ref{inequ} for details.

The core of Theorem \ref{fotetel} is a general result saying that every large
enough convex part of an expander measure keeps expanding, at least on small
enough subsets. A Borel subset $Y\in\mathcal{B}$ is $\lambda$\emph{-invariant}%
, if $e_{\lambda}(Y)=0$. When $Y$ is $\lambda$-invariant\emph{ }and $\mu
(Y)>0$, we can naturally restrict $\lambda$ to $\mathrm{Aut}$($Y,\frac{1}%
{\mu(Y)}\mu$). We say that $Y$ is $\lambda$\emph{-ergodic}, if $Y$ is
$\lambda$-invariant and every $\lambda$-invariant Borel subset $Z\subseteq Y$
satisfies $\mu(Z)=0$ or $\mu(Y\diagdown Z)=0$.

In the following we do not assume $\lambda$ and $\lambda_{i}$ to be symmetric.

\begin{lemma}
\label{mertek}Let $\lambda$ be a Borel probability measure on $\mathrm{Aut}%
(X,\mu)$, let $\kappa>\rho^{+}(\lambda)$ and let us decompose
\[
\lambda=\kappa\lambda_{1}+(1-\kappa)\lambda_{2}%
\]
where the $\lambda_{i}$ are Borel probability measures on $\mathrm{Aut}%
(X,\mu)$. Then every $\lambda_{1}$-invariant Borel subset $Y\in\mathcal{B}$ of
positive measure satisfies
\[
\mu(Y)\geq\frac{\kappa-\rho}{1-\rho}\text{.}%
\]
\newline Also, for every Borel subset $Y\in\mathcal{B}$ with
\[
0<\mu(Y)\leq\frac{1}{2}\frac{\kappa-\rho}{1-\rho}%
\]
we have
\[
\frac{e_{\lambda_{1}}(Y)}{\mu(Y)}\geq\frac{1}{2}\frac{\kappa-\rho}{\kappa
}\text{.}%
\]
In particular, one can decompose
\[
X=%
{\displaystyle\bigcup\limits_{i=1}^{n}}
X_{i}\text{ \ \ (}n\leq\frac{1-\rho}{\kappa-\rho}\text{)}%
\]
where the $X_{i}$ are $\lambda_{1}$-ergodic Borel subsets of positive measure,
and $\lambda_{1}$ on $X_{i}$ has spectral gap ($1\leq i\leq n$).
\end{lemma}

It is natural to ask whether the explicit bounds on $e_{\lambda_{1}}$ lead to
explicit bounds on $\rho^{+}(X_{i},\mathcal{B},\mu,\lambda_{1})$ in Theorem
\ref{mertek}. Even in the case when $\lambda$ is finitely supported and $X$
stays $\lambda_{1}$-ergodic, the answer is negative in general as follows from
the work of the author and Elek \cite{abelek}. However, when $X$ is
homogeneous enough, one can indeed get such explicit bounds. The nicest case
is when $X=G$ is a compact topological group and $\mu$ is the normalized Haar
measure. In this case, $G$ acts on $X$ by p.m.p. maps from both the left and
the right, and the two actions commute. The following theorem concentrates on
the case when $G$ is connected and $\kappa\geq2\rho^{+}(\lambda)$.

\begin{corollary}
\label{application}Let $G$ be a compact, connected topological group with
normalized Haar measure $\mu$. Let $\lambda$ be a Borel probability measure on
$G$ and let $\kappa\geq2\rho^{+}(\lambda)>0$. Let us decompose
\[
\lambda=\kappa\lambda_{1}+(1-\kappa)\lambda_{2}%
\]
where the $\lambda_{i}$ are Borel probability measures on $G$. Then $G$ is
$\lambda_{1}$-ergodic and we have
\[
\rho^{+}(\lambda_{1})<1-\frac{1}{512}\left(  \frac{\rho^{+}(\lambda)}{\log
_{2}(2/\rho^{+}(\lambda))}\right)  ^{2}\text{.}%
\]

\end{corollary}

In Corollary \ref{application2} we also state a version that estimates the
spectral radius $\rho(\lambda_{1})$ in terms of $\rho(\lambda)$. These
estimates have been used in the recent paper of Lindenstrauss and Varju
\cite{linvar}. Note that when $G$ is not connected, in particular, when it is
a profinite group, one can still get an explicit, but significantly weaker
estimate on the spectral gap on the ergodic components of $\lambda_{1}$, using
the machinery introduced by the author and Elek in \cite{abelek}. \bigskip

Theorem \ref{fotetel} translates to the following result in terms of profinite
actions and property $(\tau)$. For a finite $d$-regular graph $G$ let
$\rho^{+}(G)$ denote the largest nontrivial eigenvalue of the Markov (random
walk)\ operator on $G$.

\begin{theorem}
\label{profinite}Let $\Gamma$ be a group generated by the finite symmetric set
$S$. Let $(\Gamma_{n})$ be a chain of finite index subgroups in $\Gamma$ and
let
\[
\rho^{+}=\sup\rho^{+}(\mathrm{Sch}(\Gamma/\Gamma_{n},S))\text{.}%
\]
where $\mathrm{Sch}(\Gamma/\Gamma_{n},S)$ is the Schreier graph of the coset
action on $\Gamma/\Gamma_{n}$. If $H\leq\Gamma$ is a subgroup such that
\[
\rho(\Gamma,H,\lambda_{S})>\rho^{+}%
\]
then there exists a finitely generated subgroup $H^{\prime}\leq H$ such that
$(H^{\prime}\cap\Gamma_{n})$ has property ($\tau$) in $H^{\prime}$. If the
$\Gamma_{n}$ are normal in $\Gamma$, then there exists $M>0$ such that
$\left\vert \Gamma:H^{\prime}\Gamma_{n}\right\vert $ $<M$. 
\end{theorem}

A special case of Theorem \ref{profinite} is Shalom's theorem \cite{shalom}:
he assumes $\rho^{+}=\rho(T_{S})$ and that $H$ is a nontrivial normal subgroup
of $\Gamma$. Here $T_{S}$ is the $\left\vert S\right\vert $-regular tree.
Indeed, in this case, by Kesten's theorem \cite{kesten}, $\Gamma$ must be a
free product of cyclic groups and so $H$ is nonamenable, which, again using
Kesten's theorem, implies $\rho(\Gamma,H,\lambda_{S})>\rho^{+}$.

A major motivation of Shalom's result was to find infinite index subgroups of
arithmetic groups that have property $(\tau)$ with respect to its congruence
subgroups. An early provocative question in this direction was the so-called
$1$-$2$-$3$ problem of Lubotzky. In the last decade, this arithmetic direction
experienced an enormous activity, starting in the breakthrough works of
Helfgott \cite{helfgott} and Bourgain-Gamburd \cite{bourgamb} and continued in
a series of deep papers. See \cite{bourgamb}, \cite{bourgvar},
\cite{golsefvarju}, \cite{pyberszabo}, \cite{bregreentao}\ for the latest
developments. In particular, proprty ($\tau$) is now known to hold for a large
class of arithmetic lattices in semisimple Lie groups acting on their
congruence completion, assuming $H$ is Zariski dense. We do not expect that
Theorems \ref{fotetel} and \ref{profinite} will say much new in this
direction, because it seems rather nontrivial to effectively estimate the
spectral gap of a profinite action. An advantage of our result is that it
substitutes the arithmetic language, namely Zariski density and congruence
subgroups with a simple spectral condition, in the spirit of Gamburd
\cite{gamburd}, but in the discrete setting. \bigskip

\noindent\textbf{Ramanujan actions. }Let $\lambda$ be a symmetric probability
measure on the countable group $\Gamma$ and let $\varphi$ be a p.m.p. action
of $\Gamma$. We call the triple $(\Gamma,\lambda,\varphi)$ \emph{Ramanujan},
if $\rho(\lambda)=\rho^{+}(\varphi(\lambda))$. Theorem \ref{fotetel} implies
that for a Ramanujan action, every subgroup $H$ where $\rho(\Gamma
,H,\lambda)>\rho(\Gamma,1,\lambda)$ acts with finitely many ergodic components
and spectral gap on each components. These actions seem to be tight in many
other senses. For instance, a recent theorem of the author, Glasner and Virag
\cite{irsramanujan} plus an even more recent result of Bader, Duchesne and
Lecureux \cite{badule} implies that for such actions, the stabilizer of a
$\mu$-random $x\in X$ in $\Gamma$ lies in the amenable radical of $\Gamma$
a.s. Indeed, as we show in Proposition \ref{locglob}, $\rho^{+}(\varphi
(\lambda))\geq\rho(\Gamma,H,\lambda)$ where $H$ is a typical stabilizer. Being
Ramanujan then implies that $\rho(\Gamma,H,\lambda)=\rho(\lambda)$, which, by
\cite{irsramanujan} implies that $H$ is amenable a.s. Now \cite{badule}
implies that every stabilizer of a p.m.p. action that is amenable a.s. lies in
the amenable radical.

Straightforward examples for Ramanujan actions are nontrivial Bernoulli
actions (or, in another language, i.i.d. processes) of $\Gamma$
\cite{kechristsankov}. Note that for a Bernoulli action of $\Gamma$, every
infinite subgroup $H$ acts ergodically and the restricted action is also a
Bernoulli, in particular, it has spectral gap if and only if $H$ is
nonamenable. So for Bernoulli actions Theorem \ref{fotetel} does not give much
new. Another, quite non-trivial examples come from the Lubotzky-Philips-Sarnak
construction \cite{lps} that produces Ramanujan profinite actions for free
groups, for suitable ranks and the standard generating set. Recently,
Backhausz, Szegedy and Virag \cite{backszevir} analyzed the behavior of local
algorithms using the notion of a Ramanujan graphing. The connection is that
when $\lambda=\lambda_{S}$ for a symmetric generating set, the triple
$(\Gamma,\lambda,\varphi)$ is encoded in a graphing. \bigskip

\noindent\textbf{Acknowledgements. }We thank P\'{e}ter Varj\'{u} for helpful
discussions, in particular, for suggesting a more streamlined proof for
Corollary \ref{application2}. \bigskip

The paper is organized as follows. In Section \ref{prelimsection} we define
the basic notions and prove some general lemmas on expansion and spectral gap.
In particular, we prove Lemma \ref{schmidt}, a weaker substitute for Schmidt's
lemma \cite{schmidt} for probability measures on $\mathrm{Aut}(X,\mu)$.
Section \ref{pmpsection} contains the proofs of Lemma \ref{mertek} and Theorem
\ref{fotetel}. In Section \ref{profinitesection} we prove Corollary
\ref{application} and Corollary \ref{application2} and discuss profinite and
Ramanujan actions.

\section{Preliminaries \label{prelimsection}}

In this section we define the basic notions and state some general lemmas.

We start with a proposition that is well-known in Riemannian geometry and
finite graph theory under the name Cheeger inequalities. For the nontrivial
part we use the exposition by Lyons and Nazarov \cite{lynaz}.

\begin{proposition}
\label{inequ}Let $(X,\mathcal{B},\mu)$ be a standard Borel probability space
and let $\lambda$ be a Borel probability measure on $\mathrm{Aut}(X,\mu)$.
Then we have
\[
1-\mathrm{h}(\lambda)\leq\rho^{+}(\lambda)\leq1-\frac{\mathrm{h}(\lambda)^{2}%
}{8}\text{.}%
\]

\end{proposition}

\noindent\textbf{Proof. }For $Y\in\mathcal{B}$ let $f_{Y}=\chi_{Y}-y\chi_{X}$
and let $y=\mu(Y)$. Then $f\in L_{0}^{2}(X,\mu)$ and
\[
e_{\lambda}(Y)=\left\langle \chi_{Y}M_{\lambda},\chi_{X}-\chi_{Y}\right\rangle
=y(1-y)-\left\langle f_{Y}M_{\lambda},f_{Y}\right\rangle \text{.}%
\]
This gives
\[
\mathrm{h}(\lambda)\geq\frac{e_{\lambda}(Y)}{y(1-y)}=1-\frac{\left\langle
f_{Y}M_{\lambda},f_{Y}\right\rangle }{\left\langle f_{Y},f_{Y}\right\rangle
}\geq1-\rho^{+}(\lambda)\text{.}%
\]
For the other inequality, see e.g. \cite[Theorem 3.1]{lynaz}, where this is
proved in the case when $\lambda$ is the uniform measure on a finite subset
$S$ of $\mathrm{Aut}(X,\mu)$. The proof therein goes through without
difficulty by formally changing
\[
\frac{1}{\left\vert S\right\vert }\sum\limits_{s\in S}\left(  .\right)  \text{
\ to }\int\limits_{\mathrm{Aut}(X,\mu)}\left(  .\right)  d\lambda
\]
everywhere. $\square$\bigskip

Curiously, we did not find a non-probabilistic proof for the following.

\begin{proposition}
\label{nulla}Let $(X,\mathcal{B},\mu)$ be a standard Borel probability space
and let $\lambda$ be a Borel probability measure on $\mathrm{Aut}(X,\mu)$.
Then $\rho^{+}(\lambda)\geq0$.
\end{proposition}

\noindent\textbf{Proof. }We can assume that $X$ is the unit circle and $\mu$
is the normalized Lebesque measure. Let $d$ be the normalized distance on $X$,
so that the length of $X$ is $1$.\ Let $B(x,r)$ be the ball of radius $r$
around $x$.

For a parameter $r>0$ let us define the random vector
\[
f_{r}=\chi_{B(x,r)}-\chi_{B(y,r)}%
\]
where $x,y$ are independent $\mu$-random elements of $X$. Then $f_{r}\in
L_{0}^{2}(X,\mu)$, $\left\langle f_{r},f_{r}\right\rangle \leq4r$ and
\[
\mathcal{P}(\left\langle f_{r},f_{r}\right\rangle <4r)\leq4r\text{.}%
\]

Let $\gamma\in\mathrm{Aut}(X,\mu)$ be fixed. The expected measure of the
intersection $B(x,r)\cap B(y,r)g$ can be computed by fixing $y$ and only using
that $B(y,r)g$ has $\mu$-measure $2r$. This gives
\[
\mathbb{E}(\mu(B(x,r)\cap B(y,r)\gamma))=4r^{2}\text{.}%
\]
This implies
\[
\left\langle f_{r}\gamma,f_{r}\right\rangle \geq-\left(  \mu(B(x,r)\cap
B(y,r)\gamma+\mu(B(x,r)\gamma\cap B(y,r)\right)
\]
which yields
\[
\mathbb{E}(\left\langle f_{r}\gamma,f_{r}\right\rangle )\geq-8r^{2}\text{.}%
\]

Let $g_{r}=f_{r}$ conditioned on $d(x,y)\geq2r$. This implies $\left\langle
g_{r},g_{r}\right\rangle =4r.$ The probability of this event is $(1-4r)$, so
using $\left\vert \left\langle f_{r}\gamma,f_{r}\right\rangle \right\vert
\leq4r$, we have%

\[
\mathbb{E}\left(  \left\langle g_{r}\gamma,g_{r}\right\rangle \right)
\geq-\frac{1}{1-4r}(8r^{2}+16r^{2})=-\frac{24r^{2}}{1-4r}\text{.}%
\]

Let $\gamma\in\mathrm{Aut}(X,\mu)$ be $\lambda$-random, independently of
$g_{r}$. We get
\[
\mathbb{E}_{g_{r}}\left(  \left\langle g_{r}M_{\lambda},g_{r}\right\rangle
\right)  =\mathbb{E}_{g_{r}}\mathbb{E}_{\gamma}\left(  \left\langle
g_{r}\gamma,g_{r}\right\rangle \right)  =\mathbb{E}_{\gamma}\mathbb{E}_{g_{r}%
}\left(  \left\langle g_{r}\gamma,g_{r}\right\rangle \right)  \geq
-\frac{24r^{2}}{1-4r}\text{.}%
\]

In particular, for every $r>0$ there exists $g\in L_{0}^{2}(X,\mu)$ such that
\[
\frac{\left\langle gM_{\lambda},g\right\rangle }{\left\langle g,g\right\rangle
}\geq\frac{-24r}{1-4r}\text{.}%
\]
Letting $r$ tend to zero implies $\rho^{+}(\lambda)\geq0$. $\square$\bigskip

We will use the following easy lemma multiple times. We omit the proof.

\begin{lemma}
\label{metszet}Let $(X,\mathcal{B},\mu)$ be a standard Borel probability space
and let $\lambda$ be a Borel probability measure on $\mathrm{Aut}(X,\mu)$.
Then for all $Y,Z\in\mathcal{B}$ we have
\[
e_{\lambda}(Y\cap Z)\leq e_{\lambda}(Y)+e_{\lambda}(Z)\text{.}%
\]
Similarly, $e_{\lambda}(Y\backslash Z)\leq e_{\lambda}(Y)+e_{\lambda}(Z)$.
\end{lemma}

In the following lemma we prove that for an ergodic measure, if all small sets
expand, then the measure is an expander. For p.m.p. actions of countable
groups, this is proved e.g. in \cite{abelek}, using Schmidt's lemma
\cite{schmidt}. Since we could not find a version of this lemma for Borel
measures on $\mathrm{Aut}(X,\mu)$, we give a direct proof here that is of
combinatorial nature.

\begin{lemma}
\label{schmidt}Let $(X,\mathcal{B},\mu)$ be a standard Borel probability
space. Let $\lambda$ be an ergodic Borel probability measure on $\mathrm{Aut}%
(X,\mu)$. Assume that there exists $c,c^{\prime}>0$ such that for every
$Y\in\mathcal{B}$ with $0<\mu(Y)<c$ we have
\[
\frac{e_{\lambda}(Y)}{\mu(Y)(1-\mu(Y))}>c^{\prime}\text{.}%
\]
Then $\rho^{+}(\lambda)<1$.
\end{lemma}

\noindent\textbf{Proof. }Let us define the function $F:[0,1]\rightarrow
\lbrack0,1]$ as
\[
F(y)=\inf\left\{  e_{\lambda}(Y)\mid Y\in\mathcal{B}\text{, }\mu(Y)=y\right\}
\text{.}%
\]
Then $F(0)=F(1)=0$ and $F$ is symmetric to $1/2$. For all $Y,Z\in\mathcal{B}$
with $Y\cap Z=\varnothing$ we have
\begin{equation}
\left\vert e_{\lambda}(Y\cup Z)-e_{\lambda}(Y)\right\vert \leq e_{\lambda
}(Z)\leq\mu(Z) \tag{Sub}%
\end{equation}
which implies that
\[
\left\vert F(y+z)-F(y)\right\vert \leq z\text{ (}y\in\lbrack0,1],z\in
\lbrack0,1-y]\text{).}%
\]
In particular, $F$ is continuous on $[0,1]$. By the assumption of the lemma,
we have
\[
F(y)>c^{\prime}y>0\text{ \ \ (}0<y<c\text{).}%
\]
Let
\[
r=\min\left\{  y\in(0,1]\mid F(y)=0\right\}  \text{ and }\varepsilon_{0}%
=\max\left\{  F(y)\mid0\leq y\leq r\right\}  \text{.}%
\]
For $0\leq\varepsilon<\varepsilon_{0}$ let
\[
g(\varepsilon)=\min\left\{  z\mid F(y)\geq\varepsilon\text{ for all }%
y\in\lbrack z,r-z]\right\}  \text{.}%
\]
Since $F(0)=F(1)=0$, $\lim_{\varepsilon\rightarrow0}g(\varepsilon)=0$. Let
$\varepsilon_{1}>0$ such that for all $0<\varepsilon\leq\varepsilon_{1}$ we
have $g(\varepsilon)\leq r/3$.

Let $Y,Z\in\mathcal{B}$ with $\mu(Y)=\mu(Z)=r$ such that
\[
e_{\lambda}(Y),e_{\lambda}(Z)\leq\varepsilon/2\leq\varepsilon_{1}/2.
\]
Using Lemma \ref{metszet}, we have
\[
\varepsilon\geq e_{\lambda}(Y)+e_{\lambda}(Z)\geq e_{\lambda}(Y\cap Z)\geq
F(\mu(Y\cap Z))
\]
which, by the definition of $g$, implies that
\[
\mu(Y\cap Z)\leq g(\varepsilon)\text{ or }\mu(Y\cap Z)\geq r-g(\varepsilon
)\text{.}%
\]
The assumption $g(\varepsilon)\leq r/3$ ensures that the above two
possibilities are mutually exclusive.

Let $\varepsilon_{1}>\varepsilon_{2}>\ldots$ be a positive sequence converging
to $0$. Since $F(r)=0$, for all $n\geq1$ there exists $Y_{n}\in\mathcal{B}$
such that $\mu(Y_{n})=r$ and $e_{\lambda}(Y_{n})\leq\varepsilon_{n}/2$. By the
above, for every pair of integers $n<m$ exactly one of the following holds:
\[
\mu(Y_{n}\cap Y_{m})\leq g(\varepsilon_{n})\text{ or }\mu(Y_{n}\cap Y_{m})\geq
r-g(\varepsilon_{n})\text{.}%
\]
We call $Y_{n}$ and $Y_{m}$ \emph{overlapping}, if $\mu(Y_{n}\cap Y_{m})\geq
r-g(\varepsilon_{n})$.

We claim that there exists an infinite subsequence of $(Y_{n})$ such that all
pairs in the subsequence are overlapping. Assume not. Then the graph defined
on $\{Y_{n}\}$ by the overlapping relation does not admit an infinite complete
subgraph, so by the infinite Ramsey theorem, it admits an infinite empty
subgraph. Let $m>2/r$ and choose $Y_{n_{1}},\ldots,Y_{n_{m}}$ from this
subgraph such that $g(\varepsilon_{n_{i}})<r/(2(m-1))$ ($1\leq i\leq m$). Then
the sets
\[
Y_{n_{i}}\backslash\bigcup\limits_{j\neq i}Y_{n_{j}}%
\]
are mutually disjoint of measure at least $r/2$, thus the sum of their
measures is at least $mr/2>1$, a contradiction. The claim holds. By passing to
this subsequence we can assume that $(Y_{n})$ is totally overlapping.

Let us endow the set of measurable subsets of $(X,\mathcal{B},\mu)$ modulo
nullsets with the usual metric
\[
d(Y,Z)=\mu(Y\backslash Z\cup Z\backslash Y)\text{.}%
\]
By the above, for all $n<m$ we have
\[
d(Y_{n},Y_{m})\leq2g(\varepsilon_{n})\text{.}%
\]
That is, in the metric $d$, $Y_{n}$ forms a Cauchy sequence. Since $d$ defines
a complete metric space, there exists $Y\in\mathcal{B}$ such that
$d(Y,Y_{n})\rightarrow0$. However, by (Sub), $e_{\lambda}$ continuous with
respect to $d$, so we have
\[
\mu(Y)=r\text{ and }e_{\lambda}(Y)=\lim_{n\rightarrow\infty}e_{\lambda}%
(Y_{n})\leq\lim_{n\rightarrow\infty}\varepsilon_{n}/2=0\text{.}%
\]
By the ergodicity of $\lambda$, this implies that $r=1$.

In particular, $F(y)>0$ for $c\leq y\leq1-c$ and so $\mathrm{h}(\lambda)>0$
and by Lemma \ref{inequ}, we have $\rho^{+}(\lambda)<1$. $\square$\bigskip

In the proof of Theorem \ref{fotetel} we use the known trick of making
$\lambda$ (and the associated random walk) lazy, by averaging it with
$\lambda_{e}$, the Dirac measure on the identity element. This gives us the
following advantages.

\begin{lemma}
\label{lazy}Let $\lambda$ be a symmetric Borel probability measure on
$\mathrm{Aut}(X,\mu)$ and let
\[
\lambda^{\prime}=\frac{1}{2}\lambda+\frac{1}{2}\lambda_{e}\text{.}%
\]
Then
\[
\frac{1}{2}\rho^{+}(\lambda)+\frac{1}{2}=\rho^{+}(\lambda^{\prime})=\left\Vert
M_{\lambda^{\prime}}\right\Vert =\left\Vert M_{\lambda^{\prime}}%
^{n}\right\Vert ^{1/n}\text{ (}n\geq1\text{)}%
\]
and
\[
\mathrm{h}(\lambda^{\prime})=\frac{1}{2}\mathrm{h}(\lambda)\text{.}%
\]

\end{lemma}

\noindent\textbf{Proof. }The spectrum of $M_{\lambda^{\prime}}$ lies in
$[0,1]$ and so $\rho^{+}(\lambda^{\prime})=\left\Vert M_{\lambda^{\prime}%
}\right\Vert =\left\Vert M_{\lambda^{\prime}}^{n}\right\Vert ^{1/n}$. By
definition, we have
\[
e_{\lambda^{\prime}}(Y)=\frac{1}{2}e_{\lambda}(Y)
\]
for all $Y\in\mathcal{B}$ which implies $\mathrm{h}(\lambda^{\prime})=\frac
{1}{2}\mathrm{h}(\lambda)$. $\square$ \bigskip

The first part of the following lemma can be found e.g. in \cite{abnik}. We
thank P. Varju for pointing out the second part.

\begin{lemma}
\label{lefedes}Let $G$ be a compact topological group with normalized Haar
measure $\mu$ and let $A,B\subseteq G$ be measurable subsets of positive
measure. Let $g$ be a $\mu$-random element of $G$. Then the expected value%
\[
\text{\textsc{E}}(\mu(Ag\cap B))=\mu(A)\mu(B)\text{. }%
\]
When $G$ is connected, for every $k\geq2$, there exists $g_{1},\ldots,g_{k}\in
G$ such that
\[
\mu(\bigcap\limits_{i=1}^{k}Ag_{i})=\mu(A)^{k}\text{.}%
\]

\end{lemma}

\noindent\textbf{Proof.} Let
\[
U=\left\{  (a,g)\in G\times G\mid a\in A\text{, }ag\in B\right\}  .
\]
Then $U$ is measurable in $G\times G$ and using Fubini's theorem, we get
\[
\mu(A)\mu(B)=\int_{a\in A}\mu(a^{-1}B)=\mu^{2}(U)=\int_{g\in G}\mu(Ag\cap
B)=E(\mu(Ag\cap B))\,\text{. }%
\]

Let $g_{1},\ldots,g_{k}\in G$ be independent $\mu$-random elements and let
$Y=\cap_{i=1}^{k}Ag_{i}$. By induction on $k$, we have
\[
\text{\textsc{E}}(\mu(Y))=\mu(A)^{k}\text{.}%
\]
By the connectedness of $G$, the set of possible values for $\mu(Y)$ is a
connected subset of $\mathbb{R}$. Hence, it must contain $\mu(A)^{k}$.
$\square$

\section{Measure preserving actions \label{pmpsection}}

This section contains the results on p.m.p. actions, in particular, we prove
Lemma \ref{mertek} and Theorem \ref{fotetel}.

We start by showing that for a p.m.p. action, the local spectral radius is
less than or equal to the global one.

\begin{proposition}
\label{locglob}Let $\Gamma$ be a countable group and let $\lambda$ be a
symmetric probability measure on $\Gamma$. Let\textrm{ }$\varphi$ be a p.m.p.
action of $\Gamma$. Then
\[
\rho^{+}(\varphi(\lambda))\geq\rho(\lambda)\text{.}%
\]

\end{proposition}

\noindent\textbf{Proof. }Both $\rho$ and $\rho^{+}(\varphi(.))$ are continuous
for the $l_{1}$ distance on the space of probability measures on $\Gamma$, so
we can assume that $\lambda$ is supported on the finite symmetric set $S$. For
$x\in X$ let $H_{x}=\mathrm{Stab}_{\Gamma}(x)$ be the stabilizer of $x$ in
$\Gamma$ and let $x\Gamma$ be the orbit of $x$ under $\Gamma$. We have%

\[
\rho(\lambda)=\lim_{n\rightarrow\infty}\sqrt[2n]{p_{e,e,2n}}\leq
\lim_{n\rightarrow\infty}\sqrt[2n]{p_{e,H_{x},2n}}=\rho(\Gamma,H_{x},\lambda)
\]
In particular, we get that for $\mu$-almost every $x\in X$, the norm of
$M_{\lambda}$ acting on $l^{2}(x\Gamma)$ is at least $\rho(\lambda)$. (We will
not use this, but by the ergodicity of $\varphi(\lambda)$, this norm is
independent of $x$).

For $B\subseteq\Gamma$ let $xB=\{xb\mid b\in B\}$. Let $\varepsilon>0$. Then
for every $x\in X$ there exists a minimal $n(x)\in\mathbb{N}$ and a function
$f_{x}:xS^{n(x)}\rightarrow\mathbb{Q}$ such that $\left\langle f_{x}%
,f_{x}\right\rangle =1$,
\[
\sum_{y\in S^{n(x)}}f_{x}(y)=0\text{ \ \ and \ \ }\left\langle f_{x}%
M_{\lambda},f_{x}\right\rangle >\rho(\lambda)-\varepsilon\text{. }%
\]
Since $n(x)$ is a measurable function of $x$, we can choose $f_{x}$ to be a
measurable function of $x$, say, by listing all the $f_{x}$-es with zero sum
and norm $1$ and taking the first one satisfying the inequality. We get that
there exists a Borel subset $Y\in\mathcal{B}$ with $\mu(Y)>0$ and $n>0$ and
such that for all $x\in Y$ we have $n(x)=n$. By a standard argument (see e.g.
\cite{kechrismiller}), by passing to a subset of positive measure, we can also
assume that
\begin{equation}
y_{1}S^{n+1}\cap y_{2}S^{n+1}=\varnothing\text{ \ (}y_{1},y_{2}\in Y\text{,
}y_{1}\neq y_{2}\text{).} \tag{Dis}%
\end{equation}
Let $F:X\rightarrow\mathbb{Q}$ be defined by
\[
F(z)=\left\{
\begin{array}
[c]{cc}%
f_{x}(z) & z\in yS^{n}\text{, }y\in Y\\
0 & \text{otherwise}%
\end{array}
\right.  .
\]
By the above, $F$ is well-defined, $F\in L_{0}^{2}(X,\mu)$ and
\[
\left\langle F,F\right\rangle =\int\limits_{y\in Y}\left\langle f_{y}%
,f_{y}\right\rangle d\mu=\mu(Y)\text{.}%
\]
Also, by (Dis), $M_{\lambda}$ acts separately on the $yS^{n+1}$ ($y\in Y$),
thus
\[
\left\langle FM_{\lambda},F\right\rangle =\int\limits_{y\in Y}\left\langle
f_{y}M_{\lambda},f_{y}\right\rangle d\mu>(\rho(\lambda)-\varepsilon
)\mu(Y)\text{.}%
\]
This implies $\rho^{+}(\varphi(\lambda))\geq\rho(\lambda)-\varepsilon$.
$\square$ \bigskip

We now prove Lemma \ref{mertek} saying that a big enough convex part of an
expander measure keeps being expander on its ergodic components. \bigskip

\noindent\textbf{Proof of Lemma \ref{mertek}.} Let $\rho=\rho^{+}(\lambda)$.
Fix $Y\in\mathcal{B}$, let $y=\mu(Y)$ and let $e_{\lambda_{1}}=e_{\lambda_{1}%
}(Y)$.

Let $f:X\rightarrow\mathbb{R}$ be defined as
\[
f(x)=\left\{
\begin{array}
[c]{cc}%
1-y & x\in Y\\
-y & x\notin Y
\end{array}
\right.  .
\]
Then $f\in L_{0}^{2}(X,\mu)$ and $\left\langle f,f\right\rangle =y(1-y)$.

For $a\in\mathrm{Aut}(X,\mu)$ let
\[
\delta_{a}=\delta_{a}(Y)=\mu(Ya\diagdown Y)\text{.}%
\]
Then we have
\[
\left\langle fa,f\right\rangle =(y-\delta_{a})(1-y)^{2}-2\delta_{a}%
y(1-y)+(1-y-\delta_{a})y^{2}=y(1-y)-\delta_{a}\text{.}%
\]
and using $\delta_{a}\leq y$ this gives
\[
y(1-y)=\left\langle f,f\right\rangle \geq\left\langle fa,f\right\rangle
=y(1-y)-\delta_{a}\geq-y^{2}\text{.}%
\]

Using the decomposition of $\lambda$, we have
\begin{align*}%
{\displaystyle\int\limits_{a\in\mathrm{Aut}(X,\mu)}}
\left\langle fa,f\right\rangle d\lambda(a)  &  =\kappa%
{\displaystyle\int\limits_{a\in\mathrm{Aut}(X,\mu)}}
\left\langle fa,f\right\rangle d\lambda_{1}(a)+(1-\kappa)%
{\displaystyle\int\limits_{a\in\mathrm{Aut}(X,\mu)}}
\left\langle fa,f\right\rangle d\lambda_{2}(a)\\
&  \geq\kappa%
{\displaystyle\int\limits_{a\in\mathrm{Aut}(X,\mu)}}
\left\langle fa,f\right\rangle d\lambda_{1}(a)-(1-\kappa)y^{2}\\
&  =-(1-\kappa)y^{2}+\kappa y(1-y)-\kappa%
{\displaystyle\int\limits_{a\in\mathrm{Aut}(X,\mu)}}
\delta_{a}d\lambda_{1}(a)\\
&  =y(\kappa-y)-\kappa e_{\lambda_{1}}\text{.}%
\end{align*}
By the definition of $\rho$ we have
\[%
{\displaystyle\int\limits_{a\in\mathrm{Aut}(X,\mu)}}
\left\langle fa,f\right\rangle d\lambda(a)\leq\rho\left\langle
f,f\right\rangle =\rho y(1-y)\text{.}%
\]
Putting together the two inequalities we get
\[
y(\kappa-y)-\kappa e_{\lambda_{1}}\leq\rho y(1-y)
\]
which yields
\begin{equation}
\frac{e_{\lambda_{1}}}{y}\geq\frac{\kappa-\rho}{\kappa}-y(\frac{1-\rho}%
{\kappa})\text{.} \tag{Exp}%
\end{equation}

In particular, if
\[
0<y<\frac{\kappa-\rho}{1-\rho}%
\]
then $e_{\lambda_{1}}>0$, so $Y$ can not be $\lambda$-invariant. This implies
that we can decompose
\[
X=%
{\displaystyle\bigcup\limits_{i=1}^{n}}
X_{i}\text{ \ \ (}n\leq\frac{1-\rho}{\kappa-\rho}\text{)}%
\]
where the $X_{i}$ are $\lambda_{1}$-ergodic Borel subsets of positive measure.

Using (Exp) again, for
\[
0<y\leq\frac{1}{2}\frac{\kappa-\rho}{1-\rho}%
\]
we get
\[
\frac{e_{\lambda_{1}}}{y}\geq\frac{1}{2}\frac{\kappa-\rho}{\kappa}\text{.}%
\]
In particular, there exists $c,c^{\prime}>0$ such that for each $X_{i}$
($1\leq i\leq n$) and every $Y\subseteq X_{i}$ with $0<\mu(Y)\leq c$ we have
$e_{\lambda_{1}}/y\geq c\prime$. By Lemma \ref{schmidt} this implies that
$\lambda_{1}$ has spectral gap on the $X_{i}$. The lemma holds. $\square
$\bigskip

We are ready to prove the main theorem of the paper. \bigskip

\noindent\textbf{Proof of Theorem \ref{fotetel}.} Let
\[
\lambda^{\prime}=\frac{1}{2}\lambda_{S}+\frac{1}{2}\lambda_{e}\text{.}%
\]
Then $\varphi(\lambda^{\prime})=\frac{1}{2}\varphi(\lambda_{S})+\frac{1}%
{2}\lambda_{e}$, by\ Lemma \ref{lazy} we have
\begin{equation}
\rho(\Gamma,H,\lambda^{\prime})-\rho^{+}(\varphi(\lambda^{\prime}))=\frac
{1}{2}\left(  \rho(\Gamma,H,\lambda^{\prime})-\rho^{+}(\varphi(\lambda
^{\prime}))\right)  >0 \tag{A}%
\end{equation}
Recall that
\[
\rho(\Gamma,H,\lambda^{\prime})=\lim_{n\rightarrow\infty}\sqrt[2n]{p_{e,H,2n}}%
\]
where
\[
p_{e,H,n}=\sum\limits_{h\in H}p_{e,h,n}=\sum\limits_{h\in H}\left\langle
\chi_{e}M_{\lambda^{\prime}}^{n},\chi_{h}\right\rangle
\]
is the probability that the $\lambda^{\prime}$-random walk of length $n$
starting at $e$ ends at $H$. By (A) there exists an even integer $n$ such
that
\[
p_{e,H,n}>\rho^{+}(\varphi(\lambda^{\prime}))^{n}\text{.}%
\]
Fix this $n$. Let $\lambda$ be the $n$-fold convolution of $\lambda^{\prime}$.
That is, $\lambda(\{g\})=p_{e,g,n}$ ($g\in\Gamma$). Let $\kappa=\lambda(H)$.
By Lemma \ref{lazy} we have
\[
\kappa=p_{e,H,n}>\rho^{+}(\varphi(\lambda^{\prime}))^{n}=\rho^{+}%
(\varphi(\lambda))\text{.}%
\]
Since $M_{\lambda}$ is self-adjoint, we have
\[
\rho^{+}(\lambda)=\rho^{+}(\lambda^{\prime})^{n}=\left\Vert M_{\lambda
^{\prime}}\right\Vert ^{n}\text{.}%
\]

Let the measures $\lambda_{1}$ and $\lambda_{2}$ on $\Gamma$ be defined by
\[
\lambda_{1}(g)=\left\{
\begin{array}
[c]{cc}%
\frac{1}{\kappa}\lambda(g) & g\in H\\
0 & g\notin H
\end{array}
\right.  \text{ and }\lambda_{2}(g)=\left\{
\begin{array}
[c]{cc}%
0 & g\in H\\
\frac{1}{1-\kappa}\lambda(g) & g\notin H
\end{array}
\right.  \text{.}%
\]
That is, we decompose $\lambda$ according to $\Gamma=H\cup(\Gamma\backslash
H)$. Then the $\lambda_{i}$ are symmetric probability measures on $\Gamma$ and
we have
\[
\lambda=\kappa\lambda_{1}+(1-\kappa)\lambda_{2}%
\]
which implies
\[
\varphi(\lambda)=\kappa\varphi(\lambda_{1})+(1-\kappa)\varphi(\lambda
_{2})\text{.}%
\]

Let $\rho=\rho^{+}(\varphi(\lambda))$. Applying Lemma \ref{mertek} on this
convex sum, we get the decomposition
\[
X=%
{\displaystyle\bigcup\limits_{i=1}^{n}}
X_{i}\text{ \ \ (}n\leq\frac{1-\rho}{\kappa-\rho}\text{)}%
\]
where the $X_{i}$ are $\varphi(\lambda_{1})$-ergodic Borel subsets of positive
measure, and the restriction of $\varphi(\lambda_{1})$ on $X_{i}$ has spectral
gap ($1\leq i\leq n$).

Since $\lambda^{\prime}$ is supported on $S\cup\{e\}$, the support $T$ of
$\lambda_{1}$ is contained in $H\cap(S\cup\{e\})^{n}$, in particular, it is
finite. Let $H^{\prime}\leq H$ be the subgroup generated by $T$. Then
$H^{\prime}$ acts on $X_{i}$ with spectral gap ($1\leq i\leq n$). $\square$
\bigskip

\noindent\textbf{Remark. }It is easy to see that finitely many bad eigenvalues
will not disturb Theorem \ref{fotetel}. More precisely, we can define the
\emph{essential top of the spectrum of }$\lambda$ as the infimum of $\rho
^{+}(\lambda)$ acting on the orto-complements of finite dimensional
$M_{\lambda}$-invariant subspaces of $L_{2}(X)$. Theorem \ref{fotetel} then
works with the essential top of the spectrum as input. \bigskip

\section{Homogeneous and profinite actions \label{profinitesection}}

In this section we prove Corollary \ref{application} and some versions of it
using norm instead of $\rho^{+}$. Then we translate Theorem \ref{fotetel} to
the profinite setting to obtain Theorem \ref{profinite}. \bigskip

\noindent\textbf{Proof of Corollary \ref{application}.} Let $X=G$ with the
standard Borel sets and let $\mu$ be the normalized Haar measure on $G$. Then
$G$ acts on itself from the left by $\mu$-preserving maps, so $\lambda$ gives
a Borel measure on $\mathrm{Aut}(X,\mu)$. Let $\rho=\rho^{+}(\lambda)>0$.

Applying Lemma \ref{mertek} we can decompose
\[
G=%
{\displaystyle\bigcup\limits_{i=1}^{n}}
X_{i}\text{ \ \ (}n\leq\frac{1-\rho}{\kappa-\rho}\text{)}%
\]
where the $X_{i}$ are $\lambda_{1}$-ergodic Borel subsets of positive measure,
and $\lambda_{1}$ on $X_{i}$ has spectral gap ($1\leq i\leq n$). Since the
right and left $G$-actions commute, for all $g\in G$ and $1\leq i\leq n$,
$X_{i}g$ is also $\lambda_{1}$-ergodic. By Lemma \ref{lefedes} there exists
$g\in G$ such that $\mu(X_{1}\cap X_{1}g)=\mu(X_{1})^{2}$. Since $\mu
(X_{1}\cap X_{1}g)$ is also $\lambda_{1}$-ergodic, we have $\mu(X_{1})=1$ and
so $G$ is $\lambda_{1}$-ergodic.

Again using Lemma \ref{mertek} and $\kappa\geq2\rho>0$ we get that for every
Borel subset $Y\in\mathcal{B}$ with
\begin{equation}
0<\mu(Y)\leq\frac{1}{2}\rho\leq\frac{1}{2}\frac{\kappa-\rho}{1-\rho}
\tag{Small}%
\end{equation}
we have
\begin{equation}
e_{\lambda_{1}}(Y)\geq\frac{1}{2}(1-\frac{\rho}{\kappa})\mu(Y)\geq\frac{1}%
{4}\mu(Y)\text{.} \tag{Exp}%
\end{equation}

Let $Z\in\mathcal{B}$ with $0<\mu(Z)\leq1/2$. Then there exists $k$ such that
\[
\frac{1}{2}\mu(Z)\rho\leq\mu(Z)^{k}\leq\frac{1}{2}\rho\text{.}%
\]
By Lemma \ref{lefedes} there exists $g_{1},\ldots,g_{k}\in G$ such that for
$Y=\cap_{i=1}^{k}Zg_{i}$ we have $\mu(Y)=\mu(Z)^{k}$.

Using Lemma \ref{metszet} and (Exp) we get
\[
ke_{\lambda_{1}}(Z)\geq e_{\lambda_{1}}(Y)\geq\frac{1}{4}\mu(Y)=\frac{1}{4}%
\mu(Z)^{k}\geq\frac{1}{8}\rho\mu(Z)
\]
which gives
\[
\frac{e_{\lambda_{1}}(Z)}{\mu(Z)(1-\mu(Z))}\geq\frac{\rho}{8k(1-\mu(Z))}%
>\frac{\rho}{8k}\geq\frac{\rho}{8\log_{2}(2/\rho)}%
\]
Since $Z$ was arbitrary and $e_{\lambda_{1}}(Z)=e_{\lambda_{1}}(Z^{c})$, this
yields
\[
h(\lambda_{1})\geq\frac{\rho}{8\log_{2}(2/\rho)}%
\]
which, by Lemma \ref{inequ}, implies
\[
\rho^{+}(\lambda_{1})\leq1-\frac{\mathrm{h}(\lambda_{1})^{2}}{8}\leq1-\frac
{1}{512}\left(  \frac{\rho^{+}(\lambda)}{\log_{2}(2/\rho^{+}(\lambda
))}\right)  ^{2}\text{.}%
\]
The Corollary holds. $\square$ \bigskip

With a more careful analysis one can certainly shave off the constant $512$.
It is not clear whether the $\log$ term can be omitted here. The square comes
from playing back and forth between $\rho^{+}$ and $h$, and thus using both
Cheeger inequalities. \bigskip

Now we prove a version that uses $\kappa$ in the estimate but in turn allows
$\rho^{+}(\lambda)=0$.

\begin{lemma}
\label{ize}Let $G$ be a compact, connected topological group with normalized
Haar measure $\mu$. Let $\lambda$ be a Borel probability measure on $G$ and
let $\kappa>2\rho^{+}(\lambda)$. Let us decompose
\[
\lambda=\kappa\lambda_{1}+(1-\kappa)\lambda_{2}%
\]
where the $\lambda_{i}$ are Borel probability measures on $G$. Then $G$ is
$\lambda_{1}$-ergodic and we have
\[
\rho^{+}(\lambda_{1})<1-\frac{1}{2^{11}}\left(  \frac{\kappa}{\log
_{2}(4/\kappa)}\right)  ^{2}\text{.}%
\]

\end{lemma}

The proof is identical to the proof of Corollary \ref{application} above,
except that in (Small) we change $\mu(Y)\leq\rho/2$ to $\mu(Y)\leq\kappa/4$
and then throughout the whole proof, we use $\kappa/2$ instead of $\rho$ everywhere.

This leads to the following corollary, that uses $\rho$ instead of $\rho^{+}$
both as input and output. Note that the paper \cite{linvar} uses this form.

\begin{corollary}
\label{application2}Let $G$ be a compact, connected topological group with
normalized Haar measure $\mu$. Let $\lambda$ be a Borel probability measure on
$G$ and let $\kappa>2\rho(\lambda)$. Let us decompose
\[
\lambda=\kappa\lambda_{1}+(1-\kappa)\lambda_{2}%
\]
where the $\lambda_{i}$ are Borel probability measures on $G$. Then $G$ is
$\lambda_{1}$-ergodic and we have
\[
\rho(\lambda_{1})<1-\frac{1}{2^{11}}\left(  \frac{\kappa}{\log_{2}(4/\kappa
)}\right)  ^{2}\text{.}%
\]

\end{corollary}

\noindent\textbf{Proof.} For $g\in G$ let $\lambda g$ denote the $g$-translate
of $\lambda$. Then we have
\[
\rho^{+}(\lambda g)\leq\rho(\lambda)\text{ \ \ (}g\in G\text{)}%
\]
Using Lemma \ref{ize} on $\lambda g,\lambda_{1}g$ and $\lambda_{2}g$ we get
\[
\rho^{+}(\lambda_{1}g)<r=1-\frac{1}{2^{11}}\left(  \frac{\kappa}{\log
_{2}(4/\kappa)}\right)  ^{2}\text{.}%
\]
Let $f\in L_{0}^{2}(X,\mu)$ with $\left\langle f,f\right\rangle =1$. Then
\[
\left\langle fM_{\lambda_{1}},fg^{-1}\right\rangle =\left\langle
fM_{\lambda_{1}g},f\right\rangle <r\text{.}%
\]
This yields
\[
\left\langle fM_{\lambda_{1}},fM_{\lambda_{1}}\right\rangle =\int\left\langle
fM_{\lambda_{1}},fg^{-1}\right\rangle d\lambda_{1}(g)<r
\]
impying
\[
\rho(M_{\lambda_{1}})<r\text{.}%
\]
The corollary holds. $\square$ \bigskip

\noindent\textbf{Remark. }Let $\Gamma$ be a countable group and let $\lambda$
be a symmetric probability measure on $\Gamma$ with averaging operator
$M_{\lambda}$ acting on $l^{2}(\Gamma)$. In the paper we implicitely use that
$\rho^{+}(\lambda)=\left\Vert M_{\lambda}\right\Vert $. The fact itself is
folkore, but we found a couple of incomplete proofs in the literature, so we
include a sketch here. Let $\mu$ be the spectral measure of $M_{\lambda}$.
Then for all $k\geq0$, the $k$-th moment of $\mu$ equals $p_{e,e,k}$, the
probability of return in $k$ steps for the $\lambda$-random walk on $\Gamma$.
In particular, all the moments of $\mu$ are non-negative. This easily implies
that the top of the support of $\mu$ is at least the absolute value of the
bottom of the support, which then implies $\rho^{+}(\lambda)=\left\Vert
M_{\lambda}\right\Vert $. \bigskip

We are ready to prove Theorem \ref{profinite}.\bigskip

\noindent\textbf{Proof of Theorem \ref{profinite}. }Let $X$ denote the
boundary of the coset tree of $\Gamma$ with respect to $(\Gamma_{n})$. This is
the inverse limit of the coset spaces $\Gamma/\Gamma_{n}$ (see \cite{abnik}
for an exposition). Note that when the $\Gamma_{n}$ are normal, $X$ equals the
profinite completion of $\Gamma$ with respect to $(\Gamma_{n})$. Then $\Gamma$
acts on $X$ by a p.m.p. action. This action is always ergodic and has spectral
gap if and only if $(\Gamma_{n})$ has property ($\tau$) in $\Gamma$.

Let $\lambda$ be the uniform measure on the symmetric set $S$. Since
$L^{2}(X)$ as a $\Gamma$-space is the union of $L^{2}(\Gamma/\Gamma_{n})$, we
have
\[
\rho^{+}(\lambda)=\sup\rho^{+}(\Gamma/\Gamma_{n},\lambda)\text{ }%
\]
where $\rho^{+}(\Gamma/\Gamma_{n},\lambda)$ is the top of the spectrum of
$\lambda$ acting on $L_{0}^{2}(\Gamma/\Gamma_{n})$. Similarly, we have
$\rho(\lambda)=\sup\rho(\Gamma/\Gamma_{n},\lambda)$.

Since $\rho(\Gamma,H,\lambda)>\rho^{+}$, we can apply Theorem \ref{fotetel}
and get that there exists a finitely generated subgroup $H^{\prime}\leq H$
such that $X$ has finitely many ergodic $H^{\prime}$-components and
$H^{\prime}$ acts on each component with spectral gap. This is equivalent to
saying that there exists $M>0$ such that the action of $H^{\prime}$ on
$\Gamma/\Gamma_{n}$ has at most $M$ orbits and uniform spectral gap $(n\geq
0)$. In particular, $(H^{\prime}\cap\Gamma_{n})$ has property ($\tau$) in
$H^{\prime}$. When the $\Gamma_{n}$ are normal in $\Gamma$, the coset action
of $H^{\prime}$ on $\Gamma/\Gamma_{n}$ is fixedpoint-free and the number of
its orbits equals the index $\left\vert \Gamma:H^{\prime}\Gamma_{n}\right\vert
$. $\square$ \bigskip

When the $\Gamma_{n}$ are not normal and the finite actions of $\Gamma$ are
far from regular, we do not expect that the index stays bounded in general.
Random actions on rooted trees, in the spirit of \cite{glasner} may give
counterexamples.

\end{document}